\documentclass[reqno]{amsart}
\usepackage{amsfonts}
\usepackage[active]{srcltx}
\usepackage{amsmath}
\usepackage{amssymb}
\usepackage{amsthm}
\usepackage{mathtools}
\usepackage{shuffle}
\usepackage{bm} 
\usepackage{upgreek}
\usepackage{graphicx} 
\usepackage{mathrsfs}
\usepackage{fancyhdr}
\usepackage{amsthm}
\usepackage{cases}
\usepackage{amsmath,amssymb,bm}
\usepackage{enumitem}
\usepackage[]{showlabels}
\usepackage{xcolor}

\usepackage{tikz-cd} 

\usepackage{enumitem}
\usepackage{wasysym}
\usepackage{mathscinet}
\usepackage{verbatim}
\usepackage{listings}
\usepackage{graphicx}
\usepackage[
bookmarks=true,         
bookmarksnumbered=true, 
colorlinks=true, pdfstartview=FitV, linkcolor=blue, citecolor=blue,
urlcolor=blue]{hyperref}
\usepackage{microtype}
\theoremstyle{plain}
\newtheorem{theorem}{Theorem}[section]

\newtheorem{problem}[theorem]{Problem}

\renewenvironment{proof}[1][Proof]{\textbf{#1.} }{\ \rule{0.5em}{0.5em} \par }
\theoremstyle{remark}

\theoremstyle{definition}
\newtheorem{remark}[theorem]{Remark}

\def\RR{\mathbb{R}}\def\TT{\mathbb{T}}
\def\EE{\mathbb{E}}

\def\ZZ{\mathbb{Z}}

\def\bfi{{\bf
i}}
\def\bfj{{\bf j}}

\def\cE{{\mathcal E}}
\def\cF{{\mathcal F}}

\def\la{{\lambda}}

\def\Om{{\Omega}}

\def\la{{\lambda}}

\def\bfi{{\bf i}}
\def\n{{\vec{n}}}
\def\j{{\vec{\bfj}}}
\def\l{{\vec{l}}}
\def\i{{\vec{\bfi}}}
\def\upal{{\upalpha}}
\def\upbe{{\upbeta}}
\def\umu{{\vec{\upmu}}}

\let\Section=\section
\def\section{\setcounter{equation}{0}\Section}

\title{Nonlinear Young  integrals   and  differential systems
in   H\"older media }
\date{November 6,  2019}
\author[Agrawal]{Nishant Agrawal} 
\address{Department of Mathematical and Statistical Sciences \\
 University of Alberta at Edmonton \\
Edmonton,  Canada, T6G 2G1}
\email{nishant@ualberta.ca, yaozhong@ualberta.ca, neha@ualberta.ca}

\author[Hu]{Yaozhong Hu}
 
\author[Sharma]{Neha Sharma} 

\keywords{L\'evy process, nonlinear functional of L\'evy process, multiple integrals,
chaos expansion, product formula, exponential vector, polarization technique}
\begin{document}

\title{General Product formula of multiple 
integrals of  L\'evy process  }
\footnotetext[1]{Supported by the grant NSERC RGPIN-2018-05687.} 


\subjclass[2010]{Primary: 60H05. Secondary: 60G51, 60H30}


\begin{abstract} 
We derive a  product formula for finite many  multiple stochastic integrals of L\'evy process,  expressed  in terms of the associated  Poisson random measure. The formula is compact. The proof is short and uses the exponential vectors and polarization technique.  
\end{abstract}

\maketitle

\section{Introduction} 
Stochastic analysis  of nonlinear functionals of   L\'evy processes
(including Brownian motion and Poisson process) have been studied extensively and found many applications.  There have been already many  standard  books on this topic \cite{apple,protter, sato}.  In the analysis of Brownian nonlinear functional the Wiener-It\^o  chaos expansion to expand  a nonlinear functional of Brownian motion  into the sum of multiple Wiener-It\^o integrals is a fundamental contribution to the field. The  product formula to express the product of two  (or more) multiple integrals as linear combinations of some other  multiple integrals is one of the important tools (\cite{Shi}).   It plays an important role in stochastic analysis, e.g. Malliavin calculus
(\cite{Hu, Nualart}).  

The product formula for two multiple integrals of Brownian motion is known since the work of \cite[Section 4]{Shi} and the general product formula can be found  for instance 
  in \cite[chapter 5]{Hu}.
In this paper we   give  a general formula for  the product of $m$ multiple integrals of   the Poisson random measure associated with  
(purely jump)  L\'evy process. The formula is in a compact form  and it reduced to the Shigekawa's formula when   $m=2$ and the L\'evy process is reduced to Brownian motion.

When $m=2$ a similar  formula was obtained    in \cite{LS}, where the multiple integrals is with respect to the L\'evy process itself (ours is with respect to the  associated Poisson random measure which has a better properties). To obtain their formula in \cite{LS}
Lee and Shih  use white noise analysis framework.  In this work, we only
use the classical framework in  hope that this work 
is  accessible to a different spectrum of readers. 

The product formula for multiple Wiener-It\^o formula
of Brownian motion plays an important role in many aplications 
such as   U-statistics  \cite{major}.  We hope similar things may happen.
But we are not pursuing this goal in the current  paper. 
Our formula is for purely jump L\'evy process. It can be combined with the classical formulas \cite{Hu, major, Nualart, Shi} 
to general case. 

This paper is  organized as follows. In Section 2, we give some preliminaries  on L\'evy process, the associated Poisson random measure, multiple integrals. We also state our main result in this section. In Section 3, we give the proof of the formula.

%

%

\section{Preliminary and main results} 
Let $T>0$ be a positive number and let  $\left\{\eta(t) = \eta(t, \omega)\,, 
0\le t\le T\right\}$  be  a L\'evy process on some probability space  $(\Om, \cF, P)$  with filtration $\{\cF_t\,, 0\le t\le T\}$  satisfying the usual condition.  This means that 
 $\left\{\eta(t)  \right\} $  has independent and 
 stationary increment and  the sample path is right continuous with left limit.  Without loss of generality, we assume $\eta(0)=0$.  If the process $\eta(t)$ has all moments for any time index $t$,  then presumably, one can use the polynomials of the process to approximate any nonlinear functional 
 of  the process $\left\{\eta(t)  \,, 0\le t\le T\right\}$. However, it is more convenient to use the associated  Poisson random measure to carry out the stochastic analysis   of these nonlinear functionals. 
 
The jump of the process $\eta$ at time $t$ is defined by
  \[
   \Delta \eta(t):= \eta(t) - \eta(t-) \quad \hbox{if $ \Delta \eta(t)\not=0$}\,.
   \]
Denote  $\mathbb{R}_0 := \mathbb{R} \backslash \{0\}$ and let $\mathcal{B}(\mathbb{R}_0)$ be the Borel $\sigma$-algebra generated by the family of
all Borel subsets $U \subset \mathbb{R}$, such that $\Bar{U} \subset \mathbb{R}_0$. If $U \in \mathcal{B}(\mathbb{R}_0)$ with $\Bar{U} \subset \mathbb{R}_0$ and $t>0$, we then define the {\it Poisson random measure},  
$N: [0, T]\times \mathcal{B}(\mathbb{R}_0)\times \Om\rightarrow \RR$,  associated with $\eta$ by 
\begin{equation}
N(t, U) := \sum_{0 \leq s \leq t}\chi_U(\Delta\eta(s))\,,    
\end{equation} 
where $\chi_U$ is the indicator function of $U$.  
The associated  L\'evy measure $\nu$ of $\eta$ is defined by
\begin{equation}
\nu(U) := \mathbb{E}[N(1,U)] 
\end{equation}  
and compensated jump measure $\Tilde{N}$ is defined by 
\begin{equation}
\Tilde{N}(dt,dz) := N(dt,dz) - \nu(dz)dt\,. 
\end{equation}

The stochastic integral $  \int_{\TT} f(s, z) \tilde N(ds,dz)$ 
is well-defined for a predictable process $f(s,z)$ such that 
$\int_{ \TT} \EE |f(s, z)|^2 \nu(dz)ds<\infty$,  where 
and throughout this paper we use $\TT$ 
to represent the domain   $[0, T]\times \RR_0$  to simplify notation.  

Let 
\[
 \hat L^{2, n}  := \left(L^2(\TT , \lambda\times \nu )
\right)^{\otimes n}
\subseteq L^2\left(\TT^n, \left(\lambda\times \nu  
\right)^{  n}  \right)
\]
be the space of symmetric, deterministic real functions $f$ such that
\[
\|f \|_{\hat {L}^{2, n} }^2 =  \int_{\TT^n}f^2(t_1,z_1,\cdots,t_n,z_n)dt_1\nu(dz_1)\cdots  dt_n\nu(dz_n)  < \infty \,,
\]
where $\lambda(dt)=dt$ is the Lebesgue measure. In the above the symmetry means that 
\[
f (t_1,z_1, \cdots,t_i, z_i,\cdots, t_j,z_j, \cdots,t_n,z_n)
=f (t_1,z_1,\cdots,t_j, z_j,\cdots, t_i,z_i, \cdots,t_n,z_n)
\]
for all  $1\le i<j\le n$.   
For any  $f \in \hat{L}^{2,n}$  the multiple Wiener-It\^o integral 
\begin{equation}
I_n(f) :=  \int_{\TT^n}f(t_1,z_1,\cdots,t_n,z_n)  \Tilde{N} (dt_1,dz_1)\cdots \Tilde{N} (dt_n,dz_n)  
\end{equation}
is well-defined. The importance of the introduction of the associated Poisson measure and the multiple Wiener-It\^o integrals are in the following theorem which means any nonlinear functional $F$  of the L\'evy process $\eta$ 
can be  expanded as multiple Wiener-It\^o integrals. 


\begin{theorem}[Wiener-It\^ o chaos expansion for L\'evy process] 
Let $ \cF_T =\sigma (\eta(t), 0\le t\le T)$  be $\sigma$-algebra generated the L\'evy process $\eta$.  Let $F \in L^2(\Omega, \cF_T, P)$  be an $\mathcal{F}_T $  measurable  square integrable  random variable. Then F admits the following chaos expansion: 
\begin{equation}
F = \sum_{n=0}^{\infty} I_n(f_n)\,, 
\end{equation} 
where  $f_n \in \hat{L}^{2,n}, n=1,2, \cdots  $  and where  we  denote  $I_0(f_0):=f_0=\EE(F)$.   Moreover,   we have 
\begin{equation}
\| F \|_{L^2(P)}^2 = \sum_{n=0}^{\infty}n\,! \| f_n \|_{\hat L^{2,n}}^2\,. 
\end{equation}  
\end{theorem}
This chaos expansion theorem is one of the fundamental result
in stochastic analysis of L\'evy processes.  It has been widely studied in particular when $\eta$ is the Brownian motion (Wiener process).
We refer to    \cite{Hu}, \cite{Nualart}, 
\cite{protter} and references therein for   further reading.

To state our main result of this paper, we need some notation. 
Fix a positive integer $m\ge 2$.  
Denote 
\begin{equation}
 \Upupsilon=\Upupsilon_m=\left\{ \bfi=(i_1, \cdots, i_\upal ), \ \upal=2, \cdots, m,\ 
1\le i_1<\cdots<i_\upal \le m  \right\}  \, , \label{e.def_s} 
\end{equation} 
where $\upal=|\bfi|$ is  the length of the multi-index $\bfi$
(we shall use $\upal$, $\upbe$ to denote a natural number).    It is easy to see that the cardinality of $\Upupsilon$ is   $\kappa_m:=2^m-1-m$.
Denote  $\i=(\bfi_1, \cdots, \bfi_{\kappa_m})$, which is unordered list of the elements of $\Upupsilon $, where $\bfi_\upbe\in \Upupsilon$.  We use 
  $\l= (l_{\bfi_1}, \cdots, l_{\bfi_{\kappa_m}})$  to denote a multi-index of length $ \kappa_m$  associated with $\Upupsilon$,   where $l_{\bfi_\upal}$,  
$1\le \upal \le \kappa_m$  are nonnegative integers. $\l$ can be regarded as a function from $\Upsilon$ to $ \ZZ_+=\{0, 1, 2, \cdots\}$.   Denote
\begin{equation}
\begin{cases}
\Upomega =\left\{ \l: \l:  \Upupsilon\rightarrow \ZZ_+\right\}\qquad  \hbox{and for
any $\l,\,  \n\in \Upomega$}\,,\\ 
 \upchi(k, \l, \n)
=\sum_{ 1\le \upal \le \kappa_m} 
  \left[  l_{\bfi_\upal}
\chi _{\left\{\hbox{$\bfi_\upal$ contains $k$}\right\}}  +n_{\bfi_\upal }  \chi_{
\left\{\hbox{$\bfi_\upal$ contains $k$}\right\} }    \right]\,.   
\end{cases}  
\end{equation} 
Again the above mentioned $\chi$ refers to the indicator function. The conventional notations such as    $|\l|=l_{\bfi_1}+ \cdots+ l_{\bfi_{   \kappa_m}}$; $\l!=l_{\bfi_1}! \cdots  l_{\bfi_{ \kappa_m}}\!\!  !$ and so on 
are in use.
Notice that we use $l_{\bfi_1}$ instead of $l_1$ to emphasize that the  $l_{\bfi_1}$ corresponds to $\bfi_1$.  
For $\bfi=(i_1, \cdots, i_\upal) \,,  
   \bfj=(  j_1, \cdots,   j_\upbe ) \in \Upupsilon $,  and non negative integers $\upmu$ and $\upnu$  denote 
\begin{eqnarray}
\hat\otimes_{\bfi}^{\upmu} (f_1, \cdots, f_m)
&=&\int_{([0, T]\times \RR_0)^\upmu } f_{i_1}((s_1,z_1),  \cdots, (s_\upmu, 
z_\upmu), \cdots) \hat \otimes  
\cdots \nonumber\\
&&\qquad \hat\otimes f_{i_\upal}((s_1,z_1), \cdots, 
(s_\upmu, z_\upmu), \cdots)  ds_1d\nu(z_1)\cdots 
 \nonumber\\
&&\qquad ds_\upmu dz_\upmu \   f_1\hat\otimes \cdots \hat\otimes\hat f_{i_1}\hat\otimes\cdots \hat\otimes\hat f_{i_\upal}\cdots \hat\otimes f_m\,, \label{e.2.15}
\end{eqnarray}
and 
\begin{eqnarray}
  && V_{  \bfj}^\upnu (f_1, \cdots, f_m)
=  f_{j_1}((t_1,z_1),  \cdots, (t_\upnu  , z_\upnu  ), \cdots) \hat \otimes  
\cdots  \nonumber\\
&&\qquad   \hat\otimes f_{j_\upbe }((t_1,z_1),  \cdots, (t_\upnu , z_\upnu  ),  \cdots)    f_1\hat\otimes \cdots \hat\otimes\hat f_{j_1}\hat\otimes\cdots \hat\otimes\hat f_{j_\upbe }\cdots \hat\otimes f_m\,, 
\label{e.2.16} 
\end{eqnarray}
where $\hat\otimes$ denotes the symmetric tensor product and 
$ \hat f_{j_1}$ means that the function $f_{j_1}$ is removed  from the list.
Let us emphasize that both $\hat\otimes_{\bfi}^{\upmu}$ and  $V_{  \bfj}^\upnu$ are well-defined when the lengths of  $\bfi$ and $\bfj$ are one. However, we shall not use $\hat\otimes_{\bfi}^{\upmu}$ when $|\bfi|=1$  and when $|\bfj|=1$,
$V_{  \bfj}^\upnu (f_1, \cdots, f_m)
=    f_1\hat\otimes \cdots \hat\otimes  f_m$ (namely, the identity operator). 
For any two elements $\l=(l_{\bfi_1}, \cdots, l_{\bfi_{\kappa_m}})$ and
$\n=(\mu_{\bfj_1}, \cdots, \mu_{\bfj_{ \kappa_m}})$ in $\Upomega$, denote 
\begin{equation}\
\hat\otimes_\i^\l =\hat\otimes_{\bfi_1, \cdots, \bfi_{\kappa_m}}^{l_{\bfi_1}, \cdots, l_{\bfi_{\kappa_m}}}
=\hat\otimes_{\bfi_1}^{l_{\bfi_1}}  \cdots 
  \hat\otimes_{\bfi_{\kappa_m}}^{l_{\bfi_{\kappa_m}}}\,,\quad 
 V_\j^\n= V_{ \bfj_1, \cdots,  \bfj_{  \kappa_m}}^{  \mu_{\bfj_1}, \cdots,   n_{j_{  \kappa_m}}} 
  =V_{ \bfj_1}^{  \mu_{\bfj_1}} \hat\otimes \cdots 
  \hat\otimes 
  V_{  \bfj_{   \kappa_m}}^{  n_{ \bfj_{ \kappa_m}}}\,. \label{e.2.19} 
\end{equation} 
Now we can state the main result of the paper. 
\begin{theorem}\label{t.2.2}
Let $f_k \in  \left(L^2([0, T]\times \RR_0, dt\otimes \nu(dz) )\right)^{\hat\otimes q_k}$, $k=1, \cdots, m$.   
Then  
\begin{eqnarray} 
&&\prod_{k=1}^m  I_{q_k}(f_k)
=  \sum_{{{\l, \n\in \Upomega}\atop{{\upchi(1, \l,\n)\le q_1}\atop{\cdots} }}\atop
 {\upchi(m, \l,\n)\le q_m}}
\frac{\prod_{k=1}^m q_k!}{\prod_{\upal=1}^{\kappa_m}  l_{\bfi_\upal}! 
\prod_{\upbe=1}^{\kappa_m} \mu_{\bfj_\upbe}!   \prod_{k=1}^{ m}    (q_k-\upchi(k, \l,\n))! }\nonumber\\
&&\qquad \qquad 
I_{|q|+|\n|-|\upchi(  \l,\n)|}
(\hat\otimes_{\bfi_1, \cdots, \bfi_{\kappa_m}}^{l_{\bfi_1}, \cdots, l_{\bfi_{\kappa_m}}}  \hat\otimes  
  V_{ \bfj_1, \cdots,  \bfj_{  \kappa_m}}^{  \mu_{\bfj_1}, \cdots,   n_{j_{  \kappa_m}}}  (f_1, \cdots, f_m ))
    \,, 
      \label{e.2.main_formula}  
\end{eqnarray} 
where we recall $|q|=q_1+\cdots+q_m$ and $|\upchi(\l,\n)|=\upchi(1,\l,\n)+\cdots+\upchi(m, \l,\n)$.  
\end{theorem}
If $m=2$,  then $\kappa_m= 1$.  To shorten the notations we can write $q_1=n$, $q_2=m$, 
$f_1=f$,  $f_2=g$,   $l_{\upal_1}=l$, $n_{\upbe_1}=k$. Thus, $\upchi(1,\l,\n)= \upchi(2,\l,\n)=
l+k$  and $|q|+|\n|-|\upchi(\l,\n)|=n+m+k-2(l+k)=n+m-2l-k$. Thus for the product of two multiple integrals  the above theorem can be written as 
\begin{theorem}\label{t.2.3}
Let $f\in \left(L^2([0, T]\times \RR_0, dt\otimes \nu(dz) )\right)^{\hat\otimes     
n}$  and $g\in \left(L^2([0, T]\times \RR_0, dt\otimes \nu(dz) )\right)^{\hat\otimes  m}$. 
Then  
\begin{eqnarray}
 I_n(f_n)I_m(g_m)
 &=&     \sum_{\substack{k,l\in \ZZ_+\\
 k+l\le m\wedge n}} \frac{n!m!}{l! k! (n-k-l)!(m-k-l)!}  I_{ n+m-2l-k }\Big( f_n \hat{\otimes}_{k,l}  g_m \Big) \,, \nonumber\\
\label{e.2.7} 
\end{eqnarray}  
where  $\ZZ_+$ denotes the set of   non negative integers  and 
\begin{eqnarray}
&&f_n\hat \otimes_{k, l} 
 g_m(s_1, z_1, \cdots, s_{n+m-k-2l}, z_{n+m-k-2l})\nonumber\\
&&\qquad ={\rm symmetrization \ of}\ \
\int_{\TT^l} f_n(s_1, z_1,
\cdots, s_{n -l}, z_{n -l}, t_1,  y_1, \cdots, t_l, y_l) \nonumber\\
&&\qquad\qquad \qquad   g_m(s_1, z_1, \cdots, s_k, z_k,   s_{n-l+1}, \cdots, z_{n-l+1},
\cdots, \nonumber\\
&&\qquad\qquad \qquad  s_{n+m-k-2l}, \cdots, z_{n+m-k-2l}, t_1, z_1, \cdots,   t_l, z_l)
dt_1\nu(dz_1)\cdots dt_l\nu(dz_l)\,. \nonumber\\
\label{e.2.8} 
\end{eqnarray}
\end{theorem}
\begin{remark}
\begin{enumerate}
\item[(1)] 
When $\eta$ is the Brownian motion, the product formula \eqref{e.2.7} is known since 
\cite{Shi} (see e.g. \cite[Theorem 5.6]{Hu} for a formula of the general  form 
\eqref{e.2.main_formula})  and is given by 
\begin{eqnarray}
 I_n(f_n)I_m(g_m)
 &=&     \sum_{ l=0}^{n\wedge m}\frac{n!m!}{l!  (n-l)!(m -l)!}  I_{ n+m-2l  }\Big( f_n \hat{\otimes}_{ l}  g_m \Big) \,. \nonumber\\
\label{e.2.9} 
\end{eqnarray}
It is a ``special case" of \eqref{e.2.7} when $k=0$.
\end{enumerate}   
\end{remark}

\section{Proof of  Theorem \ref{t.2.2}  }
We shall  prove the main result (Theorem \ref{t.2.2})   of this paper.   
We shall prove this by using the polarization technique 
(see \cite[Section 5.2]{Hu}).  
First, let us find the  Wiener-It\^o chaos expansion for the {\it exponential functional} (random variable)  of the form 
\begin{eqnarray} \label{e.3.1} 
  Y(T)&=&  \cE(\rho(s,z))\nonumber\\
    &:=& \exp\left\{ \int_\TT \rho(s,z)\Tilde{N}(dz,ds) -\int_\TT \Big( e^{\rho(s,z)}-1-\rho(s,z)\Big)\nu(dz)ds    \right\} \nonumber\\
\end{eqnarray}
where   $\rho(s,z) \in \hat L^{2 }:=\hat L^{2, 1}=L^2(\TT, \nu(dz)\otimes \la(dt))$.  
An application of It\^o formula (see e.g. \cite{ protter}) yields 
\[
Y(T) = 1+ \int\limits_{0}^T\int\limits_{\mathbb{R}_0}Y(s-)\Big[ \exp{(\rho(s,z))} - 1 \Big]\Tilde{N}(ds,dz)\,.
\]
Repeatedly using this formula, we   obtain the chaos expansion of $Y(T)$ 
as follows. 
\begin{eqnarray}
  \cE(\rho(s,z)) 
 &=&  \exp\left\{ \int_\TT \rho(s,z)\Tilde{N}(dz,ds) -\int_\TT \Big( e^{\rho(s,z)}-1-\rho(s,z)\Big)\nu(dz)ds    \right\}\nonumber\\
 &= &\sum_{n=0}^{\infty}\frac{1}{n!} I_n(f_n)\,, \label{e.3.2}
\end{eqnarray}
where  the convergence is in $L^2(\Om, \cF_T, P)$   and 
\begin{equation} 
    f_n  =   f_n( s_1,z_1,\cdots,s_n,z_n) =(e^{\rho}-1)^{\hat\otimes n}=
     \prod_{i=1}^n\big(e^{\rho( s_i,z_i)} -1\big) \,.  \label{e.3.3}
\end{equation} 
%
%
We shall first make critical application of the above expansion
formula \eqref{e.3.2}-\eqref{e.3.3}.  For any   functions $p_k(s,z) \in \hat L^2$  (in what follows when we write  $k$ we always mean  $k=1, 2, \cdots, m$ and we shall omit $k=1, 2, \cdots, m$), 
   we  denote  
\begin{equation}
\rho_k(u_k,s,z)=\log(1+u_k p_k(s,z)) \,, 
\label{e.3.4} 
\end{equation}  

From \eqref{e.3.2}-\eqref{e.3.3}, we have (consider $u_k$   as   fixed 
real numbers) 
\begin{equation} 
\cE(\rho_k(u_k, s,  z))=\sum_{n=0}^{\infty} \frac1{n!} u_k^n I_n(f_{k, n} ) \,, 
 \end{equation}  
where 
\begin{equation} 
f_{k, n} = \frac{1}{u_k^n} \prod_{i=0}^n (e^{\rho_k(u_k,s_i,z_i)} - 1 ) =p_k^{\otimes n}= \prod_{i=1}^n p_k(s_i,z_i) 
\label{e.3.6}
\end{equation}  
%
%
It is clear that 
\begin{equation} 
\prod_{k=1}^m \cE(\rho_k(u_k, s,  z)) =\sum_{q_1,\cdots, q_m=0}^{\infty} \frac1{q_1!\cdots q_m!} u_1^{q_1}\cdots u_m^{q_m}  I_{q_1}(f_{1, q_1}  )\cdots  I_{q_m}(f_{m,q_m}  )
 \,,  \label{e.3.7} 
 \end{equation}   
where $f_{k, q_k}$, $k=1, \cdots, m$     are defined by  \eqref{e.3.6}.

On the other hand, from the definition  of the 
exponential functional \eqref{e.3.1},  
we have
\begin{eqnarray}
&&\prod_{k=1}^m \cE(\rho_k(u_k, s,  z)) \nonumber \\
&&\qquad  = \prod_{k=1}^m  \exp{\Big\{   \int_\TT \rho_k(u_k, s,  z)\Tilde{N}(dz,ds) -\int_\TT\Big( e^{\rho_k(u_k, s,  z)}-1-\rho_k(u_k, s,  z)\Big)\nu(dz)ds   \Big\}}  \nonumber \\ 
& &\qquad = \exp \Big\{ \int_\TT \sum_{k=1}^m \rho_k(u_k, s,  z) \Tilde{N}(dz,ds) \nonumber \\
&&\qquad\qquad -\int_\TT\Big( e^{\sum_{k=1}^m \rho_k(u_k, s,  z)} -1-\sum_{k=1}^m \rho_k(u_k, s,  z) \Big)\nu(dz)ds    \Big\}  \nonumber \\  
&&\qquad\qquad  \cdot \exp    \Big\{ \int_\TT e^{ \sum_{k=1}^m \rho_k(u_k, s,  z)}-\sum_{k=1}^m e^{\rho_k(u_k, s,  z)}  +m-1\Big)\nu(dz)ds 
 \Big\}
 \nonumber \\
&&\qquad =:A\cdot B\,,  \label{e.3.8} 
\end{eqnarray}
where $A$ and $B$ denote the above first and second exponential factors. 

The first exponential factor $A$ is an exponential functional of the form \eqref{e.3.1}. Thus, again by the chaos expansion formula 
\eqref{e.3.2}-\eqref{e.3.3},  we have  
\begin{equation}
A = \sum_{ {n}=0}^{\infty} \frac1{n!} I_{ {n}}(h_{ {n}}(u_1, \cdots, u_m)) \,,
\label{e.3.9} 
 \end{equation} 
 where
 \begin{eqnarray}
 &&h_{ {n}}(u_1, \cdots, u_m)
=   \prod_{i=0}^{  {n}}  (
 e^{\sum_{k=1}^m \rho_k(u_k,s_i,z_i) } - 1 )\,. \label{e.3.10} 
  \end{eqnarray} 
By the definition of $q_k$,   we have 
\begin{eqnarray*}
\sum_{k=1}^m \rho_k(u_k,s_i,z_i)
&=&\log\prod_{k=1}^m  (1+u_kp_k(s_i,z_i) ) \,.  
\end{eqnarray*} 
Or 
\def\Sym{{\rm Sym}}
\begin{eqnarray*}
h_{ {n}}(u_1, \cdots, u_m)
&=&\left(\left[ \prod_{k=1}^m  (1+u_kp_k  )-1\right]\right)^{\hat \otimes n}\\
    &=& \Sym_{(s_1,z_1), 
    \cdots, (s_n, z_n)} \prod_{i=1}^{n }\left[ \prod_{k=1}^m  (1+u_kp_k(s_i,z_i) )-1\right]\,,  
\end{eqnarray*}
where $\hat \otimes $ denotes the symmetric tensor product and 
$\Sym_{(s_1,z_1), 
    \cdots, (s_n, z_n)}$ denotes the symmetriization with respect to 
    $(s_1,z_1), 
    \cdots, (s_n, z_n)$.
\def\bfj{{\bf j } } 
Define
\[
S=\left\{ \bfj=(j_1, \cdots, j_\upbe ), \ \upbe=1, \cdots, m,\ 
1\le j_1<\cdots<j_\upbe \le m  \right\}\,. 
\]
The cardinality of $S$ is  $|S|=\tilde \kappa_m:=2^m-1$. 
We shall freely use the notations  introduced in Section 2. Denote also 
\begin{eqnarray*} u_\bfj=u_{j_1}\cdots u_{j_\upbe}\,,\quad 
p_\bfj(s,z)=p_{j_1}(s,z)\cdots p_{j_\upbe}(s,z)
\quad(\hbox{for $\bfj=(j_1, \cdots, j_\upbe) \in S$})\,.  
\end{eqnarray*} 
We have    
\begin{eqnarray*}
&&h_{ {n}}(u_1, \cdots, u_m) 
= \left(  \sum_{\bfj\in S} u_\bfj p_\bfj  \right)^{\hat \otimes n}
=\sum_{ |\umu|=n}\frac{|\umu|!}{\umu!}u_\j^\umu p_\j^{\hat\otimes \umu}\\
&&\qquad\qquad =  \sum  _{\mu_{\bfj_1}+\cdots+\mu_{\bfj_{\tilde\kappa_m}}=n}
\frac{n!}{
\mu_{\bfj_1}!   \cdots  \mu_{\bfj_{\tilde\kappa_m}}!}  u_{\bfj_1}^{\mu_{\bfj_1}}   \cdots u_{\bfj_{\tilde\kappa_m}}^{\mu_{\bfj_{\tilde\kappa_m}}}  p_{\bfj_1}^{\hat\otimes \mu_{\bfj_1}}
\hat\otimes\cdots \hat\otimes p_{\bfj_{\tilde\kappa_m}}^{\hat\otimes \mu_{\bfj_{\tilde\kappa_m}}}\,,  
\end{eqnarray*}
where $\umu:S\rightarrow \ZZ_+$ is a multi-index and we used  the notation 
$u_\j^\umu=   u_{\bfj_1}^{\mu_{\bfj_1}}   \cdots u_{\bfj_{\tilde\kappa_m}}^{\mu_{\bfj_{\tilde\kappa_m}}}$; and $p_\j^{\hat\otimes \umu}= p_{\bfj_1}^{\hat\otimes \mu_{\bfj_1}}
\hat\otimes\cdots \hat\otimes p_{\bfj_{\tilde\kappa_m}}^{\hat\otimes \mu_{\bfj_{\tilde\kappa_m}}}$.  
Inserting the above expression   into 
\eqref{e.3.9} we have 
\begin{eqnarray}
A 
&=& \sum_{ {n}=0}^{\infty} \sum  _{\mu_{\bfj_1}+\cdots+\mu_{\bfj_{\tilde\kappa_m}}=n}
\frac{1}{
\mu_{\bfj_1}!   \cdots  \mu_{\bfj_{\tilde\kappa_m}}!}  u_{\bfj_1}^{\mu_{\bfj_1}}   \cdots u_{\bfj_{\tilde\kappa_m}}^{\mu_{\bfj_{\tilde\kappa_m}}} 
I_n( p_{\bfj_1}^{\hat\otimes \mu_{\bfj_1}}
\hat\otimes\cdots \hat\otimes p_{\bfj_{\tilde\kappa_m}}^{\hat\otimes \mu_{\bfj_{\tilde\kappa_m}}})\nonumber\\
\label{e.3.11} 
 \end{eqnarray}
Now we consider the second exponential factor in \eqref{e.3.8}:
\begin{align*}
B&=
\exp    \Big\{ \int_\TT \Big( e^{ \sum_{k=1}^m \rho_k(u_k, s,  z)}-\sum_{k=1}^m e^{\rho_k(u_k, s,  z)}  +m-1\Big)\nu(dz)ds 
 \Big\} \nonumber \\
&=\exp \Big\{ \sum_{\bfi\in \Upupsilon } u_\bfi \int_{\TT} p_\bfi (s,z)   \nu(dz)ds  \Big\}  \,,
\end{align*}
where $\Upupsilon $ is defined by \eqref{e.def_s}  
(which is a subset of $S$ such that $|\bfj|\ge 2$).    
Thus,   
\begin{eqnarray}
B &=&
\sum_{n=0}^\infty \frac{1}{n!} \left( \sum_{\bfi\in \Upupsilon } u_\bfi \int_{\TT} p_\bfi (s,z)   \nu(dz)ds  \right)^n \nonumber  \\
&=& \sum_{n=0}^\infty \sum_{ {l_{\bfi_1}+\cdots+l_{\bfi_{\kappa_m}}=n} } \frac{1}{l _{\bfi_1}!
\cdots l_{\bfi_{\kappa_m}}!}   u_{\bfi_1}^{l_{\bfi_1}}\cdots 
u_{\bfi_{\kappa_m}}^{l_{\bfi_{\kappa_m}}} \left(\int_{\TT} p_{\bfi_1} (s,z)   \nu(dz)ds
\right)^{l_{\bfi_1} } \nonumber\\
&&\qquad\qquad \cdots \left(\int_{\TT} p_{\bfi_{\kappa_m}} (s,z)   \nu(dz)ds \right)^{l_{\bfi_{\kappa_m}} }  \,,  
\label{e.3.12} 
\end{eqnarray}
where $\l\in \Upomega$ is a multi-index.  
Combining \eqref{e.3.11}-\eqref{e.3.12}, we have 
\begin{eqnarray}
&& AB
 = \sum_{n, \tilde n=0}^\infty \sum_{ {\mu_{\bfj_1}+\cdots+\mu_{\bfj_{\tilde\kappa_m}}=n }\atop  {l_{ \bfi_1}+\cdots+l_ 
{ \bfi_{\kappa_m}} =\tilde  n}} 
\frac{1}{
\mu_{\bfj_1}!   \cdots  \mu_{\bfj_{\tilde\kappa_m}}! l_{\bfi_1}!
\cdots l_{\bfi_{\kappa_m}}!}    u_{\bfj_1}^{\mu_{\bfj_1}}   \cdots u_{\bfj_{\tilde\kappa_m}}^{\mu_{\bfj_{\tilde\kappa_m}}} \nonumber\\
&&\qquad\qquad   u_{  \bfi_1}^{l_{  \bfi_1}}\cdots 
u_{ \bfi_{\kappa_m}}^{l_{l   \bfi_{ \kappa_m}}} B_{\bfi,   \bfj, 
l_{\bfi},  \mu_{\bfj}}\,, \qquad{\rm where}\label{e.AB} 
\\
&&B_{\bfi,   \bfj, l_\bfi, \mu_\bfj}
 :=  \left(\int_{\TT} p_{  \bfi_1} (s,z)   \nu(dz)ds
\right)^{l_{\bfi_1} }  \cdots \nonumber\\
&&\qquad\qquad  \left(\int_{\TT} p_{ \bfi_{\kappa_m}} (s,z)   \nu(dz)ds \right)^{l_{\bfi_{\kappa_m}} }  
I_n( p_{\bfj_1}^{\hat\otimes \mu_{\bfj_1}}
\hat\otimes\cdots \hat\otimes p_{\bfj_{\tilde\kappa_m}}^{\hat\otimes \mu_{\bfj_{\tilde\kappa_m}}}) \,. \label{e.3.15}
\end{eqnarray}
To get an expression for 
$B_{\bfi,   \bfj, l_\bfi, \mu_\bfj}$ we use the notations 
\eqref{e.2.15}-\eqref{e.2.16} and \eqref{e.2.19}. 
Then 
\begin{equation}
B_{\bfj, \tilde \bfj, n_\bfj, \tilde n_\bfj}=I_n(\hat\otimes_{ \j }
^{  \umu }
\hat \otimes V_\i^\l
 (p_1^{\otimes n_{\bfi_1}}, \cdots, p_m^{\otimes n_m}))\,. \label{e.3.20} 
\end{equation}
To compare the coefficients of 
$u_1^{n_{\bfi_1}}\cdots u_m^{n_m}$, we need to express the right hand side of \eqref{e.AB}  as a power series of $u_1, \cdots, u_m$.  
For $k=1, \cdots , m$   denote  
\begin{equation}
 \tilde \upchi(k, \l, \umu)
=\sum_{ 1\le \upal \le \kappa_m} 
    l_{\bfi_\upal}
I_{\left\{\hbox{$\bfi_\upal$ contains $k$}\right\}}  + \sum_{ 1\le \upbe \le \tilde \kappa_m}  \mu_{\bfj_\upbe }  I_{
\left\{\hbox{$\bfj_\upbe$ contains $k$} \right\} }      
\,.\label{e.3.21} 
\end{equation}
Combining \eqref{e.3.8}, \eqref{e.AB} and \eqref{e.3.20}, we have 
\begin{eqnarray} 
&&\sum_{q_1,\cdots, q_m=0}^{\infty} \frac{u_1^{q_1}\cdots u_m^{q_m}} {q_1!\cdots q_m!}   I_{q_1}(p_1^{\otimes q_1}  )\cdots  I_{q_m}(p_m^{\otimes q_m }   )
\nonumber\\ 
& & = \sum_{n, \tilde n=0}^\infty \sum_{ {{\mu_{\bfj_1}+\cdots+\mu_{\bfj_{\tilde\kappa_m}}=n }\atop  {l_{\bfi_1}+\cdots+l_{\kappa_m=\tilde  n}}}
\atop {   \tilde \upchi(k,\l,\umu)=q_k, k=1, \dots, m} }
\frac{u_1^{q_1}\cdots u_m^{q_m} }{l_{\bfi_1}!
\cdots l_{\bfi_{\kappa_m}}! 
\mu_{\bfj_1}!   \cdots  \upmu_{\bfj_{\tilde\kappa_m}}!  }\nonumber\\
&&\qquad \qquad I_n(\hat\otimes_{\bfi_1, \cdots, \bfi_{\kappa_m} }^{l_{\bfi_1}, \cdots, l_{\bfi_{\kappa_m}}} \hat\otimes  
  V_{  \bfj_1, \cdots, \bfj_{\tilde \kappa_m}}^{\upmu_{\bfj_1}, \cdots,   \upmu_{\bfj_{\tilde  \kappa_m}}}  (p_1^{\otimes q_1}, \cdots, p_m^{\otimes q_m}))
    \,.\nonumber\\
      \label{e.3.22} 
\end{eqnarray} 
Comparing the coefficient of $u_1^{q_1}\cdots u_m^{q_m}$, we have 
\begin{eqnarray} 
&&   \prod_{k=1}^m I_{q_k}(p_k^{\otimes q_k }  )  =  \sum_{{\bfj_1,\cdots,  \bfj_{\tilde\kappa_m}\in S}\atop  { 
\bfi_1,\cdots,   \bfi_{\kappa_m}\in \Upupsilon }} 
\sum_{ 
\tilde \upchi(k,\l,\umu)=q_k, k=1, \dots, m}  
\frac{ q_1!\cdots q_m!   }{l_{\bfi_1}!
\cdots l_{\bfi_{\kappa_m}}! 
\mu_{\bfj_1}!   \cdots  \upmu_{\bfj_{\tilde\kappa_m}}!  }\nonumber\\
&&\qquad \qquad I_n(\hat\otimes_{\bfi_1, \cdots, \bfi_{\kappa_m} }^{l_{\bfi_1}, \cdots, l_{\bfi_{\kappa_m}}} \hat\otimes  
  V_{  \bfj_1, \cdots, \bfj_{\tilde \kappa_m}}^{\upmu_{\bfj_1}, \cdots,   \upmu_{\bfj_{\tilde  \kappa_m}}}  (p_1^{\otimes q_1}, \cdots, p_m^{\otimes q_m}))\,. 
      \label{e.3.22} 
\end{eqnarray}  
Notice that when $|\bfj|=  1$, namely, $\bfj=(k), k=1, \cdots, m$,  then 
$V_\bfj^\upmu(f_1, \cdots, f_m)=f_1\hat\otimes \cdots \hat\otimes f_m$.
We can separate these terms from 
 the remaining ones, which will satisfy $|\bfj|
 \ge 2$.  Thus, the remaining multi-indices $\bfj$'s 
  consists of the set $\Upupsilon$. 
 We can write a multi-index $\umu:S\rightarrow \ZZ_+$ as 
  $\umu =(n_{(1)}, \cdots, n_{(m)}, \n)$,  where $\n\in \Upupsilon$.  
We  also observe 
    $q_k=\tilde \upchi(k, \l, \umu)=n_{(k)}+\upchi(k, \l, \n)$.  
After   replacing   $\umu$  by $\n$, \eqref{e.3.22} 
gives \eqref{e.2.main_formula}. 
This proves Theorem \ref{t.2.2}   for $f_k=p_k^{\otimes q_k}$,
$k=1, \cdots, m$.  
By polarization technique
(see e.g. \cite[Section 5.2]{Hu}), we also know the identity 
\eqref{e.2.main_formula}  holds true for $f_k=p_{k,1}\otimes\cdots \otimes  p_{k,q_k}$,
  $p_{k,q_k}\in L^2([0, T]\times \RR_0, ds\times \nu(dz))$, $k=1, \cdots, m$. 
Because both sides of \eqref{e.2.main_formula}  are multi-linear 
with respect to $f_k$, we know  \eqref{e.2.main_formula}  holds true for 
\[ 
f_k=\sum_{\ell=1}^{\upnu_k} c_{k, \ell}  p_{k,1, \ell}\otimes\cdots \otimes  p_{k,q_k, \ell}\,,\qquad k=1, \cdots, m \,,
\]
where $c_{k,\ell}$ are constants, $p_{k, k', \ell}\in  L^2([0, T]\times \RR_0, ds\times \nu(dz))$,
$k=1, \cdots, m, k'=1, \cdots, q_k$ and $\ell=1, \cdots, \upnu_k$.  
Finally, the identity \eqref{e.2.main_formula}   is proved by a 
routine limiting argument.

\end{document}